\documentclass[notitlepage,11pt]{article}
\usepackage{amssymb,amsmath
}
\catcode`\@=11
\@addtoreset{equation}{section}

\catcode`\@=12
\usepackage{colortbl}%

\def\R{\mathbb{R}}

\def\f{\varphi}

\def\irn{\int\limits_{\R^n}}

\def\Dshalf{\left(-\Delta\right)^{\!\frac{s}{2}}\!}

\def\Ds{\left(-\Delta\right)^{\!s}\!}

\def\eps{\varepsilon}

\def\proof{\noindent{\textbf{Proof. }}}
\def\QED{\hfill {$\square$}\goodbreak \medskip}

\newtheorem{Theorem}{Theorem}[section]
\newtheorem{Lemma}[Theorem]{Lemma}

\newtheorem{Corollary}[Theorem]{Corollary}
\newtheorem{Remark}[Theorem]{Remark}
\newtheorem{Definition}[Theorem]{Definition}

\begin{document}

\title{\vspace{-10mm}
{Variational inequalities for the fractional Laplacian}}

\author{Roberta Musina\footnote{Dipartimento di Matematica ed Informatica, Universit\`a di Udine,
via delle Scienze, 206 -- 33100 Udine, Italy. Email: {roberta.musina@uniud.it}. 
{Partially supported by Miur-PRIN 201274FYK7\_004.}}~ ,
Alexander I. Nazarov\footnote{
St.Petersburg Department of Steklov Institute, Fontanka 27, St.Petersburg, 191023, Russia, 
and St.Petersburg State University, 
Universitetskii pr. 28, St.Petersburg, 198504, Russia. E-mail: al.il.nazarov@gmail.com.
Supported by RFBR grant 14-01-00534 and by
St.Petersburg University grant 6.38.670.2013.}
~ and
{Konijeti Sreenadh\footnote{Department of Mathematics, Indian Institute of Technology Delhi, Hauz Khas, New Delhi-110016, India.
Email:{sreenadh@maths.iitd.ac.in}}}}

\date{}

\maketitle

\begin{abstract}
\footnotesize
In this paper we study the obstacle problems for the fractional Lapalcian
of order $s\in(0,1)$ in a bounded domain $\Omega\subset\R^n$,
under mild assumptions on the data.
%
%
\end{abstract}

\normalsize

\bigskip

\section{Introduction}
Let $\Omega$ be a bounded domain in $\R^n$, $n\ge 1$. 
Given $s\in(0,1)$, a measurable function $\psi$ and a distribution $f$ on ${\Omega}$, we consider the problem
\begin{equation}
{\label{eq:problem}}
\begin{cases}
u\ge \psi&\text{in ${\Omega}$}\\
\Ds u\ge f&\text{in ${\Omega}$}\\
\Ds u= f&\text{in $\{u>\psi\}$}\\
u=  0&\text{in  $\R^n\setminus\overline\Omega$.}
\end{cases}
\end{equation}
Our interest is motivated by the noticeable paper  \cite{Sil}, where Louis E. Silvestre
investigated (\ref{eq:problem}) in case $\Omega=\R^n$, $f=0$ and
$\psi$  smooth. His results apply also to Dirichlet's problems on balls,
see \cite[Section 1.3]{Sil}.
Besides remarkable results, in \cite{Sil} the interested reader can find  stimulating
motivations for (\ref{eq:problem}), arising from mathematical finance. 
In addition, Signorini's problem, also known as {the lower dimensional obstacle problemÓ for the classical Laplacian},
can be recovered from (\ref{eq:problem}) by taking $s=\frac{1}{2}$.

Among the papers dealing with (\ref{eq:problem}) and related problems we cite also
\cite{BCRo, CaFi, CSS, GP, PP, SV.obs} and references there-in, with no attempt to 
provide a complete reference list. 

In the present paper we show that the free boundary problem
(\ref{eq:problem}) admits a solution under quite mild assumptions on the data, 
see Theorems \ref{T:regularity} and \ref{T:2} below.  However, our starting interest included
broader questions concerning  the variational inequality
\begin{equation*}\tag{$\mathcal P(\psi,f)$}
\label{eq:vi}
u\in K^s_\psi~,\qquad
\langle \Ds  u-f, v-u\rangle\ge 0\quad\forall v\in K^s_\psi~\!,
\end{equation*}
where $f\in \widetilde H^s(\Omega)'$ and
$$
K^s_\psi=\left\{v\in{\widetilde H^s(\Omega)}~|~v\ge \psi ~\textit{a.e. on ${{\Omega}}$}~\right\}.
$$
Notation and main definitions are listed at the end of this  introduction. 
We will always assume that the closed and convex set $K^s_\psi$
is not empty, also when not explicitly stated. 

Problem \ref{eq:vi} admits a unique solution $u$, that can be characterized as
the unique minimizer for
\begin{equation}
\label{eq:minimization}
\inf_{v\in K^s_\psi}~\frac{1}{2}\langle \Ds  v,v\rangle-\langle f, v\rangle~\!.
\end{equation} 

The variational inequality \ref{eq:vi} and the free boundary problem (\ref{eq:problem})
are naturally related. Any solution $u\in \widetilde H^s(\Omega)$ to (\ref{eq:problem})
coincides with the unique solution to \ref{eq:vi}, see Remark \ref{R:uniqueness}.
Conversely, if $u$ solves \ref{eq:vi} then
$\left(-\Delta\right)^{\!s}u\!-f$ is a nonnegative distribution on $\Omega$, 
compare with Theorem \ref{T:sup}. By analogy with the local case $s=1$
one can guess that $\Ds u=f$ outside the
coincidence set $\{u=\psi\}$, at least when $u$ is regular enough. This is essentially
the content of  Section 3 in \cite{Sil}, where $f=0$ and 
$\psi$ is a smooth, rapidly decreasing function on  $\Omega=\R^n$,
and of Theorems \ref{T:regularity}, \ref{T:2} below.

To study the variational inequality \ref{eq:vi} we took inspiration from the classical 
theory about the local case $s=1$. In particular, we refer to
the fundamental monograph \cite{KS} by Kinderlehrer and Stampacchia,
and to the pioneering papers \cite{BrSt, LS1, LS2, JLL, Mi, U0, U}, among others.

Standard techniques do not apply directly in the fractional case, mostly
because of the different behavior of the truncation operator $v\mapsto v^+$, $H^s(\R^n)\to H^s(\R^n)$.
Section \ref{S:truncations} is entirely devoted to this subject; we collect there some lemmata
that might have an independent interest.

We take advantage of the results in Section \ref{S:truncations} to
obtain  equivalent and useful formulations for \ref{eq:vi}, 
and to prove continuous dependence theorems
upon the data $f$ and $\psi$, see Sections \ref{S:equivalent} and \ref{S:dependance}, respectively.

Some extra difficulties arise from having settled a nonlocal problem on a
bounded domain, producing at least, but not only, the same (partially solved) technical difficulties
as for the unconstrained problem $\Ds u=f$, $u\in \widetilde H^s(\Omega)$
(see for instance \cite{Co}, \cite{RoS}, \cite{ROS2} and references
there-in, for regularity issues).

\medskip

Our main results proved in Section \ref{S:regularity1}. 
They involve the unique solution $\omega_f$ to
\begin{equation}
\label{eq:e}
\Ds\omega_f=f\quad\text{in $\Omega$}~,\qquad \omega_f\in \widetilde H^s(\Omega).
\end{equation}

\begin{Theorem}
\label{T:regularity}
Assume that $\psi$ and $f\in\widetilde H^s(\Omega)'$ satisfy the following conditions:
\begin{itemize}
\item[$A1)$] $(\psi-\omega_f)^+\in \widetilde H^s(\Omega)$; 
\item[$A2)$] $\Ds(\psi-\omega_f)^+-f$ is a  locally finite signed measure on $\Omega$;
\item[$A3)$] $(\Ds(\psi-\omega_f)^+-f)^+\in L^p_{\rm loc}(\Omega)$ for some $p\in[1,\infty]$.
\end{itemize}
Let $u\in{\widetilde H^s(\Omega)}$ be the unique solution to \ref{eq:vi}.
Then the following facts hold.
\begin{itemize}
\item[$i)$] $\Ds u-f\in L^p_{\rm loc}({\Omega})$;
\item[$ii)$] $0\le \Ds u-f\le (\Ds(\psi-\omega_f)^+-f)^+$ a.e. on $\Omega$;
\item[$iii)$] $\Ds  u=f$ a.e. on $\{u>\psi\}$.
\end{itemize}
In particular, $u$ solves the free boundary problem (\ref{eq:problem}).
\end{Theorem}

\begin{Theorem}
\label{T:2}
Assume that $\Omega$ is a bounded Lipschitz domain satisfying the exterior ball condition.
Let $\psi\in C^0(\overline\Omega)$ be a given obstacle, such that $K^s_\psi$ is not empty, 
$\psi\le 0$ on $\partial\Omega$
and $f\in L^p(\Omega)$, for some exponent $p>n/2s$.

Then the unique solution $u$ to \ref{eq:vi} is continuous on $\R^n$ and solves the free boundary problem (\ref{eq:problem}).
\end{Theorem}

Our results
plainly cover the non-homogeneous Dirichlet's free boundary problem
$$
\begin{cases}
u\ge \psi&\text{in ${\Omega}$}\\
\Ds u\ge f&\text{in ${\Omega}$}\\
\Ds u= f&\text{in $\{u>\psi\}$}\\
u=  g&\text{in  $\R^n\setminus\overline\Omega$,}
\end{cases}
$$
under appropriate assumptions on the datum  $g$.  Notice indeed that $u$ solves the related
variational inequality if and only if
$u-g$ solves $\mathcal P(\psi-g,f+\Ds g)$.

Free boundary problems for the operator $\Ds u+u$
can be considered as well, with minor modifications in the statements and in the proofs.

\bigskip
\small
\noindent
{\bf Notation}
\label{SS:notation}
The definition of the  fractional Laplacian $\Ds$
involves the Fourier transform: 
$$
\mathcal F\left[\Ds u\right]=|\xi|^{2s}\mathcal F[u]~\!,
\quad 
{\cal F}[u](\xi)=\frac{1}{(2\pi)^{n/2}}\irn e^{-i~\!\!\xi\cdot x}u(x)~\!dx~\!.
$$
{Let $\Omega\subset\R^n$ be a bounded domain. }We adopt the  standard notation
\begin{eqnarray*}
&&{H^s(\R^n)}=\{u\in {L^2(\R^n)}~|~\Dshalf u\in L^2(\R^n)~\},\\
&&~{\widetilde H^s(\Omega)}~\!=~\!\{u\in {H^s(\R^n)}~|~u\equiv 0~~\text{on}~~\R^n\setminus\overline\Omega\}.
\end{eqnarray*}
We endow $H^s(\R^n)$ and ${\widetilde H^s(\Omega)}$ with their natural Hilbertian structures. We
recall that  the norm of $u$ in ${\widetilde H^s(\Omega)}$ is given by  the $L^2(\R^n)$-norm of $\Dshalf u$.

{We do not make any assumption on $\Omega$. Thus $\partial\Omega$ might be
very irregular,  even a fractal, and  $C^\infty_0(\Omega)$ might be not dense in 
$\widetilde H^s(\Omega)$. Notice that ${\widetilde H^s(\Omega)}$ coincides with
${\widetilde H^s(\Omega')}$,   whenever $\overline\Omega=\overline\Omega'$.}

We denote by $\langle\cdot,\cdot\rangle$ the  duality product between ${\widetilde H^s(\Omega)}$
and its dual $\widetilde H^s(\Omega)'$. In particular,
$\Ds u\in\widetilde H^s(\Omega)'$ for any $u\in \widetilde H^s(\Omega)$, and
$$
\langle\Ds u, v\rangle= 
\int\limits_{\R^n}\Dshalf u\cdot \Dshalf v~\!dx
=\int\limits_{\R^n}|\xi|^{2s}\mathcal F[u]~\!\overline{\mathcal F[v]}~\!d\xi.
$$

\normalsize

\section{Truncations}
\label{S:truncations}

For measurable functions $v,w$ we put, as usual,
$$
v\vee w=\max\{v,w\}~,\quad v\wedge w=\min\{v,w\}~,\quad v^+=v\vee 0~,\quad v_-=-(v\wedge 0),
$$
so that $v=v^+-v^-$. 
It is well known that $v\vee w\in {H^s(\R^n)}$ and $v\wedge w\in {H^s(\R^n)}$
if $v,w\in {H^s(\R^n)}$.

\begin{Lemma}
\label{L:truncation}
Let $v\in {H^s(\R^n)}$. Then 
\begin{description}
\item$~~i)$\quad $\displaystyle{
\langle \Ds  v^+,v^-\rangle=\langle \Ds v^-,v^+\rangle\le 0}$~\!;
\item$~ii)$\quad $\displaystyle{
\langle \Ds  v,v^-\rangle+ \int\limits_{\mathbb R^n}|\Dshalf v^-|^2~\!dx}\le 0$~\!;
\item$iii)$\quad $\displaystyle{
\langle \Ds  v,v^+\rangle- \int\limits_{\mathbb R^n}|\Dshalf v^+|^2~\!dx}\ge0$~\!.
\end{description}
In addition, 
if $v\in {H^s(\R^n)}$ does not have constant sign, then all the above inequalities are strict.
\end{Lemma}

\proof
In \cite[Theorem 6]{MNHS},  the Caffarelli-Silvestre extension argument \cite{CaSi}
has been used to check that
$$
\int\limits_{\mathbb R^n}|\Dshalf |v||^2~\!dx<\int\limits_{\mathbb R^n}|\Dshalf v|^2~\!dx~\!,
$$
whenever $v$ changes sign. That is,
$$
\int\limits_{\mathbb R^n}|\Dshalf(v^++v^-)|^2~\!dx<\int\limits_{\mathbb R^n}|\Dshalf(v^+-v^-)|^2~\!dx~\!.
$$
The conclusion is immediate.
\QED

\begin{Remark}
\label{R:MP}
One can use $ii)$ in Lemma \ref{L:truncation} to get the 
well known weak maximum principle, that is, if $u\in \widetilde H^s(\Omega)$
and $\Ds u\ge 0$ in $\Omega$ then $u\ge 0$ in $\Omega$.
\end{Remark}

\begin{Corollary}
\label{C:continuity}
Let $v_h$ be a sequence in $H^s(\R^n)$ such that $v_h$ converges to a nonpositive function  in $H^s(\R^n)$. Then $v_h^+\to 0$ in $H^s(\R^n)$.
\end{Corollary}

\proof
Statement  $iii)$ in Lemma \ref{L:truncation} provides the estimate
\begin{equation}
\label{eq:sw}
\int\limits_{\mathbb R^n}|\Dshalf v_h^+|^2~\!dx\le \langle \Ds  v_h,v_h^+\rangle,
\end{equation}
that gives us the boundedness of the sequence 
$v_h^+$ in $H^s(\R^n)$. Since $v_h^+\to 0$ in $L^2(\R^n)$, we have
$v_h^+\to 0$ weakly in $H^s(\R^n)$. Thus 
$\langle \Ds  v_h,v_h^+\rangle\to 0$, as $\Ds v_h$ converges  in $H^s(\R^n)'$,
and the conclusion follows from (\ref{eq:sw}). 
\QED

\begin{Lemma}
\label{L:m}
Let $v\in {\widetilde H^s(\Omega)}$ and  $m>0$.
Then 
$(v+m)^-, (v-m)^+, v\wedge m\in {\widetilde H^s(\Omega)}$ and
\begin{itemize}
\item[$i)$]\quad $\displaystyle{\langle {\Ds} v,(v+m)^-\rangle+\int\limits_{\mathbb R^n}|\Dshalf (v+m)^-|^2~\!dx\le 0}$;
\item[$ii)$]\quad $\displaystyle{\langle {\Ds} v,(v-m)^+\rangle- \int\limits_{\mathbb R^n}|\Dshalf (v-m)^+|^2~\!dx
\ge 0}$;
\item[$iii)$]\quad$\displaystyle{\int\limits_{\mathbb R^n}|\Dshalf (v\wedge m)|^2~\!dx\le 
\int\limits_{\mathbb R^n}|\Dshalf v|^2~\!dx- \int\limits_{\mathbb R^n}|\Dshalf (v- m)^+|^2~\!dx}$.
\end{itemize}
\end{Lemma}

\proof
Clearly, $(v+m)^-\in L^2(\R^n)$ and $(v+m)^-\equiv 0$ outside $\Omega$. 
Fix a cutoff function $\eta\in C^\infty_0(\R^n)$, 
with $0\le\eta\le 1$, and such that
$\eta\equiv 1$ in a  ball containing $\overline\Omega$. 
Then $(v+m)^-=(v+m\eta)^-\in \widetilde H^s(\Omega)$,
as trivially $m\eta\in  H^s(\R^n)$.

 For any integer $h\ge 1$ we set
$$
\eta_h(x)=\eta\Big(\frac{x}{h}\Big),
$$
so that $\eta_h\to 1$ pointwise. A direct computation shows that
\begin{equation}
\label{eq:etah}
\Ds \eta_h(x)=h^{-2s}\Big(\Ds \eta\Big) \Big(\frac{x}{h}\Big)~\longrightarrow~0
\quad\text{in
{$L^2_{\rm loc}(\R^n)$.}}
\end{equation}
By $ii)$ in Lemma \ref{L:truncation} we have that
\begin{eqnarray*}
0&\ge&
\langle \Ds(v+m\eta_h),(v+m)^-\rangle+\int\limits_{\mathbb R^n}|\Dshalf (v+m)^-|^2~\!dx\\
&=&\langle \Ds v,(v+m)^-\rangle+\int\limits_{\mathbb R^n}|\Dshalf (v+m)^-|^2~\!dx+
m \int_{\Omega} (\Ds \eta_h)(v+m)^-~\!dx\\
&=&\langle \Ds v,(v+m)^-\rangle+\int\limits_{\mathbb R^n}|\Dshalf (v+m)^-|^2~\!dx+o(1),
\end{eqnarray*}
by (\ref{eq:etah}) and since $(v+m)^-$ has compact support in ${\Omega}$.
Claim $i)$ is proved.
To check $ii)$ notice that $(v-m)^+=((-v)+m)^-$ and then use $i)$ with 
$(-v)$ instead of $v$.

It remains to prove $iii)$. Notice that $v\wedge m=v-(v-m)^+$. Hence
$v\wedge m\in \widetilde H^s(\Omega)$. Using  $ii)$ we get
\begin{eqnarray*}
\|\Dshalf (v\wedge m)\|^2
&=&\|\Dshalf v\|^2
-2
\langle\Ds v,(v-m)^+\rangle
+ \|\Dshalf (v-m)^+\|^2\\
&\le& 
\|\Dshalf v\|^2-
\|\Dshalf (v-m)^+\|^2~\!.
\end{eqnarray*}
The proof is complete.
\QED

\section{Equivalent formulations}
\label{S:equivalent}
We start this section by introducing a crucial notion.
 \begin{Definition}
A function ${\cal U}\in\widetilde H^s(\Omega)$ is a supersolution for $\Ds v=f$ if
$$
\langle \Ds{\cal U}-f,\f\rangle\ge  0\qquad\text{for any $\f\in \widetilde H^s(\Omega)$, $\f\ge 0$.}
$$
\end{Definition}
The above definition  extends the usually adopted one in the local case
$s=1$, see  \cite[Definition 6.3]{KS}. A different definition of supersolution is used
in \cite{Sil} for $f=0$. We refer to \cite[Subsection 2.10]{Sil}, 
for a stimulating discussion on this subject.

\begin{Theorem} 
\label{T:sup}
Let $u\in K^s_\psi$. The following sentences are equivalent.
\begin{itemize}
\item[$a)$] $u$ is the solution to problem \ref{eq:vi};
\item[$b)$] $u$ is the smallest supersolution for $\Ds v=f$ in the convex set $K^s_\psi$. That is, 
${\cal U}\ge u$ almost everywhere in ${\Omega}$,  for any supersolution ${\cal U}\in K^s_\psi$;
\item[$c)$]  $u$ is a supersolution for $\Ds v=f$ and 
$$\displaystyle{\langle\Ds u-f,(v-u)^-\rangle= 0}\qquad\text{for any $v\in K^s_\psi$.}
$$
\item[$d)$] $\displaystyle{\langle\Ds v-f,v-u\rangle\ge 0}$ for any $v\in K^s_\psi$.
\end{itemize}
\end{Theorem}

\proof
{$a) \Longleftrightarrow b)$}.
Assume that $u$ solves \ref{eq:vi}. Fix any 
nonnegative $\f\in \widetilde H^s({\Omega})$. Testing \ref{eq:vi} with 
$u+\f\in K^s_\psi$ one gets
$\langle \Ds  u-f,\f\rangle\ge0$, 
that proves that $u$ is a supersolution. 

Next, take any
supersolution ${\cal U}\in K^s_\psi$. 
Then $u-(u-{\cal U})^+={\cal U}\wedge u \in K^s_\psi$. Thus
$$
\langle \Ds  u-f,-(u-{\cal U})^+\rangle\ge 0.
$$
On the other hand, from $\Ds{\cal U}- f\ge 0$ we get $$\langle \Ds  {\cal U}-f,(u-{\cal U})^+\rangle\ge 0.$$
Adding the above inequalities we arrive at
$$
0\ge \langle \Ds (u-{\cal U}),(u-{\cal U})^+\rangle\ge 
\int\limits_{\mathbb R^n}|\Dshalf (u-{\cal U})^+|^2~\!dx,
$$
thanks to $iii)$ in Lemma \ref{L:truncation}.
Thus $(u-{\cal U})^+=0$ almost everywhere in ${\Omega}$, that is, $u\le {\cal U}$ and proves that 
$a)$ implies $b)$. 

Conversely, assume that $u$ satisfies $b)$ and let $\tilde u$ be the solution to \ref{eq:vi}. 
We already know that $a) \Rightarrow b)$. Thus $u$ and $\tilde u$ must coincide, because both
obey the condition of being the smallest supersolution to $\Ds v=f$ in $K^s_\psi$.
Hence, $a)$ holds.

\medskip

\noindent
$a) \Longleftrightarrow c)$.
Let $u$ be the solution to \ref{eq:vi}. We already know that $u$ is  supersolution.
Fix any function $v\in K^s_\psi$. 
Notice that 
$$u+(v-u)^-\ge u\ge\psi~,\quad u-(v-u)^-=v\wedge u\ge \psi.$$ Thus,
testing \ref{eq:vi} with  $u\pm(v-u)^-$ we get
$\displaystyle{\langle\Ds u-f,\pm(v-u)^-\rangle\ge  0}$,
that is, $c)$ holds. 

Conversely, assume that $u$ satisfies $c)$.
Let $\tilde u\in K^s_\psi$ be the solution to \ref{eq:vi}. We already proved that 
$\tilde u$ is the smallest supersolution in $K^s_\psi$. In particular, $\tilde u\le u$
and thus
$$
\langle\Ds u-f,u-\tilde u\rangle=\langle\Ds u-f,(\tilde u-u)^-\rangle= 0
$$
by the assumption  $c)$ on $u$. Since $\tilde u$ solves \ref{eq:vi},
we also get
$$
\langle\Ds \tilde u-f,u-\tilde u\rangle\ge  0~\!.
$$
Substracting, we infer $\langle\Ds (u-\tilde u),u-\tilde u\rangle\le  0$, that is, $u=\tilde u$. 
\\

\noindent
$a) \Longleftrightarrow d)$. Clearly $a)$ implies $d)$ because
\begin{multline*}
\langle\Ds v-f,v-u\rangle\\=\langle\Ds u-f,v-u\rangle+\langle\Ds (v-u),v-u\rangle\ge
\langle\Ds u-f,v-u\rangle~\!.
\end{multline*}
Now assume that $u$ satisfies $d)$ and fix any $v\in K^s_\psi$. 
From $\frac{v+u}{2}\in K^s_\psi$ and $d)$ we obtain
\begin{eqnarray*}
0&\le&2 \langle\Ds \Big(\frac{v+u}{2}\Big)-f,\frac{v+u}{2}-u\rangle=
\frac12\langle\Ds (v+u),v-u\rangle-\langle f,v-u\rangle\\
~\\
&=&\left(\frac12\langle\Ds v,v\rangle-\langle f,v\rangle\right)-
\left(\frac12\langle\Ds u,u\rangle-\langle f,u\rangle\right)~\!.
\end{eqnarray*}
Thus $u$ solves the minimization problem (\ref{eq:minimization}), that
is, $u$ solves \ref{eq:vi}.
\QED

\begin{Remark}
In the local case $s=1$, the equivalence between $a)$ and $d)$ is commonly known
as Minty's lemma, see \cite{Mi}.
\end{Remark}

\begin{Corollary}\label{compare_f}
Let $f_1, f_2\in \widetilde H^s(\Omega)'$ and let
 $u_i$ be  the solution to $\mathcal P(\psi,f_i)$, $i=1,2$.
If  $f_1\ge f_2$ in the sense of distributions, then $u_1\ge u_2$
a.e. in $\Omega$.
\end{Corollary}

\proof
The function $u_1$ is a supersolution for $\Ds v=f_2$ and $u_1\in K^s_{\psi}$. Hence
$u_1\ge u_2$, by statement $b)$ in Theorem \ref{T:sup}.
\QED
\begin{Remark}
\label{R:uniqueness}
Let  $u\in \widetilde H^s(\Omega)$ be a solution to (\ref{eq:problem}).
Then  $\Ds u-f$ can be identified with a nonnegative Radon measure on $\Omega$ having 
support in $\{u=\psi\}$. If $v\in K^s_\psi$, then $(v-u)^-$ vanishes on $\{u=\psi\}$.
Thus $\langle \Ds u-f,(v-u)^-\rangle=0$, hence $u$ solves \ref{eq:vi} by Theorem \ref{T:sup}.
\end{Remark}

\section{Continuous dependence results}
\label{S:dependance}

\begin{Theorem}
\label{T:bounded1}
Let $\psi_1,\psi_2$ be  given  obstacles, $f\in \widetilde H^s(\Omega)'$ and let
 $u_i$ be the  solution to $\mathcal P(\psi_i,f)$, $i=1,2$.
If $\psi_1-\psi_2\in L^\infty({\Omega})$, then $u_1-u_2$ is bounded, and
$$
i)~~\|(u_1-u_2)^+\|_\infty\le\|(\psi_1-\psi_2)^+\|_\infty~,\quad
ii)~~\|(u_1-u_2)^-\|_\infty\le\|(\psi_1-\psi_2)^-\|_\infty.
$$
\end{Theorem}

\proof
Put $m:=\|(\psi_1-\psi_2)^+\|_\infty$. Since $(u_2-u_1+m)^-\in \widetilde H^s(\Omega)$ by Lemma \ref{L:m},
then 
$$v_1:=u_1-(u_2-u_1+m)^-={(u_2+m)}\wedge u_1\in K^s_{\psi_1}~\!.
$$
Hence we can use $v_1$ as test function in $\mathcal P(\psi_1,f)$ to get
$$
\langle \Ds  u_1-f, -(u_2-u_1+m)^-\rangle\ge 0.
$$
On the other hand, we can test $\mathcal P(\psi_2,f)$ with {$u_2+(u_2-u_1+m)^-\in K^s_{\psi_2}$}. Hence
$$\langle \Ds  u_2-f, (u_2-u_1+m)^-\rangle\ge 0. $$
Adding and taking  $i)$ of Lemma \ref{L:m}  into account, we arrive at
$$
- \int\limits_{\mathbb R^n}|\Dshalf (u_2-u_1+m)^-|^2~\!dx
\ge
\langle \Ds  (u_2-u_1),(u_2-u_1+m)^-\rangle \ge 0.
$$
Hence, $(u_2-u_1+m)^-=0$. 
We have proved that $(u_1-u_2)^+\le m$ a.e. in $\Omega$, hence $i)$ holds. 
Inequality $ii)$ can be proved in the same way.
\QED

\begin{Corollary}
\label{C:infty}
Let $\psi\in  L^\infty({\Omega})$ and $f\in L^p({\Omega})$, with $p\in(1,\infty)$, $p>n/2s$.
Let $u\in{\widetilde H^s(\Omega)}$ 
be the unique solution to \ref{eq:vi}.
Then $u\in L^\infty({\Omega})$ and
\begin{equation}
\label{eq:Linf}
\psi\vee\omega_f \le u\le \|\psi^+\|_\infty+c\|f^+\|_p\quad\text{a.e. in $\Omega$,}
\end{equation}
where $\omega_f$ solves (\ref{eq:e}) 
and $c$ depends only on $n,s,p$ and  ${\Omega}$.
In particular, if $f=0$ then 
$$
\psi^+\le u\le \|\psi^+\|_\infty~\!.
$$
\end{Corollary}

\proof
First of all, notice that $f\in \widetilde H^s(\Omega)'$ by Sobolev embedding theorem.
Since $u$ is supersolution of (\ref{eq:e}), the first inequality in (\ref{eq:Linf}) follows
by the maximum principle in Remark \ref{R:MP}.

Denote by $\omega_{f^+}$ the unique solution to (\ref{eq:e}) with $f$ replaced by $f^+$. 
If $n>2s$ we use convolution to define
$$
U=c_1 |x|^{2s-n}* (f^+\cdot\chi_{\Omega}).
$$
For proper choice of the constant $c_1$, $U$ solves $\Ds U=f^+\cdot\chi_{\Omega}$ in $\mathbb R^n$.
Convolution estimates give $U\le c\|f^+\|_p$ on $\R^n$. 
{By the maximum principle,} $\omega_{f^+}\le  U$ on $\Omega$, hence $\omega_{f^+}\le
c\|f^+\|_p$. 
For $n=1\le 2s$ this inequality also holds, see, e.g., 
\cite[Remark 1.5]{RoS}.

Now let $u_1$ be the unique solution of $\mathcal P(\psi,f^+)$. Then $u_1\ge u$ by Corollary \ref{compare_f}.
Finally, we can consider $\omega_{f^+}$ as the solution of the problem $\mathcal P(\omega_{f^+},f^+)$.
Theorem \ref{T:bounded1} gives
$$
u\le (u_1-\omega_{f^+})^++\omega_{f^+}\le \|(\psi-\omega_{f^+})^+\|_\infty+\omega_{f^+},
$$
and the last inequality in (\ref{eq:Linf}) follows.
\QED

Roughly speaking, Theorem \ref{T:bounded1} concerns the continuity of 
$L^\infty\ni\psi\mapsto u\in L^\infty$.
The next  result gives the continuity of the arrow $L^\infty\ni\psi\mapsto u\in \widetilde H^s(\Omega)$.

\begin{Theorem}
\label{T:Linfty}
Let $\psi_h\in L^\infty(\Omega)$ be  a sequence of obstacles and let $f\in \widetilde H^s(\Omega)'$ be given.
Assume that there exists $v_0\in \widetilde H^s(\Omega)$, such that $v_0\ge \psi_h$ for any $h$.

Denote by  $u_h$ the solution to the obstacle problem $\mathcal P(\psi_h,f)$.
If $\psi_h\to \psi$ in $L^\infty(\Omega)$, then 
$u_h\to u$ in $\widetilde H^s(\Omega)$, where $u$ is the solution to the limiting problem
\ref{eq:vi}.
\end{Theorem}

\proof
Let $u$ be the solution to \ref{eq:vi}. We already know from Theorem \ref{T:bounded1}
that $\|u-u_h\|_\infty\le\|\psi-\psi_h\|_\infty$. Hence, in particular, $u_h\to u$ a.e. in $\Omega$.
Now, test $\mathcal P(\psi_h,f)$ with $v_0$ to obtain that
$$
\langle \Ds  u_h, u_h\rangle\le \langle \Ds  u_h- f, v_0\rangle+\langle f, u_h\rangle.
$$
Hence, the sequence $u_h$ is bounded in $\widetilde H^s(\Omega)$. Therefore,
$u_h\to u$ weakly in $\widetilde H^s(\Omega)$. To prove that $u_h\to u$ in the $\widetilde H^s(\Omega)$ norm we only
need to show that
$$
\limsup_{h\to \infty} \|\Dshalf u_h\|_2\le \|\Dshalf u\|_2.
$$
For any $\eps>0$ we introduce the function
$$
v_{\eps}=u+(v_0-u)\wedge {{\eps}}.
$$
Since $\psi_h\to \psi$ in $L^\infty(\Omega)$, we have $v_{\eps}\ge \psi_h$ for $h$ large enough.
Using $v_{\eps}$ as test function in $\mathcal P(\psi_h,f)$ we get
$$
\langle \Ds  u_h-f, u+(v_0-u)\wedge{{\eps}}-u_h\rangle\ge 0,
$$
and hence
$$
\|\Dshalf u_h\|_2^2=
\langle \Ds  u_h, u_h\rangle\le
\langle \Ds  u_h-f, u+(v_0-u)\wedge{{\eps}}\rangle+\langle f,u_h\rangle.
$$
Letting $h\to \infty$ we infer
\begin{eqnarray}
\nonumber
\limsup_{h\to \infty} \|\Dshalf u_h\|^2_2&\le&
\langle \Ds  u-f, u+(v_0-u)\wedge{{\eps}}\rangle+\langle f,u\rangle\\
&=& \|\Dshalf u\|^2_2+
\langle \Ds  u-f, (v_0-u)\wedge{{\eps}}\rangle.
\label{eq:ictp}
\end{eqnarray}
Now we let ${\eps}\to 0$. Clearly
$(v_0-u)\wedge{\eps}\to -(v_0-u)^-$ in $L^2(\Omega)$.
In addition, the functions $(v_0-u)\wedge{\eps}$ are uniformly bounded
in $\widetilde H^s(\Omega)$ by $iii)$ in Lemma \ref{L:m}. Thus 
$(v_0-u)\wedge{\eps}\to -(v_0-u)^-$ weakly
in $\widetilde H^s(\Omega)$.
Thus, 
from (\ref{eq:ictp}) we get
$$
\limsup_{h\to \infty} \|\Dshalf u_h\|^2_2 \le \|\Dshalf u\|^2_2-
\langle \Ds  u-f, (v_0-u)^-\rangle=
 \|\Dshalf u\|^2_2
$$
since $u$ solves \ref{eq:vi}, and therefore it satisfies condition $c)$ in Theorem
\ref{T:sup}. Thus $u_h\to u$ in $\widetilde H^s(\Omega)$.
\QED

Next we deal with the continuity of the arrow $H^s\ni \psi\mapsto u\in \widetilde H^s$.

 \begin{Theorem}
 \label{T:Hs2}
Let $\psi_h\in {H^s(\R^n)}$ be  a sequence of obstacles such that $\psi_h^+\in \widetilde H^s(\Omega)$,
and let $f_h$ be a sequence in 
$\widetilde H^s(\Omega)'$.
Assume that
$$
\psi_h\to \psi\quad\text{in $H^s(\R^n)$, and}\quad f_h\to f\quad\text{in $H^s(\Omega)'$}.
$$
Denote by  $u_h$ the solution to the obstacle problem
$\mathcal P(\psi_h,f_h)$. 
Then $u_h\to u$ in ${\widetilde H^s(\Omega)}$, where $u$ is the solution to the limiting obstacle problem \ref{eq:vi}.
\end{Theorem}

\proof
We can assume that $f_h, f=0$. If not, replace the obstacles $\psi_h$ and $\psi$ with
$\psi_h-\omega_{f_h}$ and $\psi-\omega_f$, respectively, see (\ref{eq:e}). 

Let $u_h$ solve $\mathcal P(\psi_h,0)$ and let
$u$ be the solution to the limiting problem $\mathcal P(\psi,0)$. Recall that $u$ is the unique 
minimizer for
\begin{equation}
\label{eq:minimization0}
\inf_{v\in K^s_\psi}~\langle \Ds  v,v\rangle~\!.
\end{equation}
Since $u\vee \psi_h=u+(\psi_h-u)^+$ and $\psi_h-u\to \psi-u\le 0$, then
\begin{equation}
\label{eq:pale2}
u\vee \psi_h\to u\quad\text{in $\widetilde H^s(\Omega)$}
\end{equation}
by Corollary \ref{C:continuity}. Moreover, $u\vee \psi_h\in K^s_{\psi_h}$ and thus from $\mathcal P(\psi_h,0)$
we infer
\begin{equation}
\label{eq:cork}
\langle \Ds  u_h, u_h\rangle\le \langle \Ds  u_h, u\vee\psi_h\rangle.
\end{equation}
Inequality (\ref{eq:cork}) guarantees the boundedness of the sequence $u_h$ in
$\widetilde H^s(\Omega)$. Hence we can assume that $u_h\to\tilde u$ weakly 
in $\widetilde H^s(\Omega)$. 
Since $\psi_h\to \psi$ and $u_h\to \tilde u$
a.e. in $\Omega$, clearly $\tilde u\in K^s_\psi$.

 Next, by weak lower semicontinuity, (\ref{eq:cork}) and (\ref{eq:pale2}) we get
\begin{equation}
\label{eq:sc}
\langle \Ds  \tilde u, \tilde u\rangle\le \liminf_{h\to\infty} \langle \Ds  u_h, u_h\rangle\le  \limsup_{h\to\infty} \langle \Ds  u_h, u_h\rangle\le \langle \Ds  \tilde u, u\rangle.
\end{equation}
Thus
$$
\|\Dshalf \tilde u\|^2_2\le \|\Dshalf \tilde u\|_2 \|\Dshalf u\|_2.
$$
Hence, $\tilde u =u$, as the minimization problem (\ref{eq:minimization0})
admits a unique solution, and
(\ref{eq:sc}) implies $\|\Dshalf u_h\|_2\to \|\Dshalf u\|_2$. Hence
$u_h\to u$  strongly
in $\widetilde H^s(\Omega)$.
\QED

\section{Proof of the main results}
\label{S:regularity1}

We start with a preliminary theorem of independent interest, that gives distributional bounds
on $\Ds u-f$ under mild assumptions on the data.

\begin{Theorem}
\label{T:measure}
Let $\psi$ and $f\in\widetilde H^s(\Omega)'$ satisfying 
assumptions $A1)$ and $A2)$ in Theorem \ref{T:regularity}.
Let $u\in{\widetilde H^s(\Omega)}$ be the unique solution to \ref{eq:vi}.
Then 
$$
0\le \Ds u-f\le (\Ds(\psi-\omega_f)^+-f)^+\quad\text{in the distributional sense on $\Omega$.}
$$
\end{Theorem}

\proof
The main tool 
was inspired by the penalty method by Lewy-Stampacchia \cite{LS1} and
already used for instance in \cite{SV.obs} under
smoothness assumptions on the data and on the solution.

\medskip

In order to simplify notations we start the proof with some remarks. 
First, we can assume that $f=0$, as we did in the proof of Theorem \ref{T:Hs2}. Thus
$\Ds u\ge 0$ and $u\ge\psi$, that imply $u\ge \psi^+$, use the maximum principle in Remark \ref{R:MP}. 
Clearly  $u$ is the smallest supersolution to $\Ds v=0$ in $K^s_{\psi^+}$, and hence it solves
the obstacle problem $\mathcal P(\psi^+,0)$. In conclusion, it suffices to prove Theorem \ref{T:measure}
in case $f=0$ and $\psi\ge 0$ in $\R^n$. Our aim
is to show that
\begin{equation}
\label{eq:ineq}
0\le \Ds u\le (\Ds\psi)^+\quad\text{in the distributional sense on $\Omega$,}
\end{equation}
for $\psi\in \widetilde H^s(\Omega)$, $\psi\ge 0$, such that $\Ds\psi$
is a measure on $\Omega$.

The  proof of (\ref{eq:ineq}) will be achieved in few steps.

\medskip

\noindent
{\bf Step 1} {\em Assume ${\Ds\psi}\in L^p(\Omega)$ for any large exponent $p>1$. Then (\ref{eq:ineq}) holds.}

\medskip
We take $p\ge \frac{2n}{n+2s}$, that is needed only if $n>2s$. Then
$\widetilde H^s(\Omega)\hookrightarrow L^{p'}(\Omega)$ and
$L^{p}(\Omega)\subset \widetilde H^s(\Omega)'$ by Sobolev embeddings. In
particular $({\Ds\psi})^+\in \widetilde H^s(\Omega)'$. 

Take a function $\theta_\eps\in C^\infty(\R)$ such that $0\le \theta_\eps\le 1$, and
$$
\theta_\eps(t)=1~~\text{for $t\le 0$,}\quad
\theta_\eps(t)=0~~\text{for $t\ge \eps$.}
$$
By standard variational methods we have that there exists a unique $u_\eps\in \widetilde H^s(\Omega)$
that weakly solves
$$
\Ds u_\eps=\theta_\eps(u_\eps-\psi)~\!({\Ds\psi})^+\quad\text{in $\Omega$.}
$$
We claim that
$$
u\le u_\eps\le u+\eps\quad\text{a.e. in $\Omega$.}
$$
By $iii)$ in Lemma \ref{L:truncation}  we can estimate
\begin{eqnarray*}
\|\Dshalf (\psi-u_\eps)^+\|_2^2&\le&
\langle \Ds(\psi-u_\eps),(\psi-u_\eps)^+\rangle\\
&\le&
\int\limits_\Omega ({\Ds\psi})^+(1-\theta_\eps(u_\eps-\psi))(\psi-u_\eps)^+~\!dx=0~\!.
\end{eqnarray*}
Hence, $u_\eps\ge \psi$. Since $\Ds u_\eps\ge 0$, then $u_\eps\ge u$
by $b)$ in Theorem \ref{T:sup}. Next, we use $iii)$ in Lemma \ref{L:m}
and $\Ds u\ge 0$
to estimate
\begin{eqnarray*}
\|\Dshalf (u_\eps-u-\eps)^+\|_2^2&\le&
\langle \Ds(u_\eps-u),(u_\eps-u-\eps)^+\rangle\\
&\le&
\int\limits_\Omega ({\Ds\psi})^+~\!\theta_\eps(u_\eps-\psi)~\!(u_\eps-u-\eps)^+~\!dx=0~\!.
\end{eqnarray*}
Thus $u_\eps\le u+\eps$, and the claim is proved. In particular, we have that
$\|u_\eps-u\|_\infty\to 0$ as $\eps\to 0$. Therefore, for any nonnegative test
function $\f\in C^\infty_0(\Omega)$ we have that
\begin{eqnarray*}
\langle \Ds u,\f\rangle&=&\int\limits_\Omega u\Ds\f~\!dx=
\int\limits_\Omega u_\eps\Ds\f ~\!dx+o(1)\\
&=&\langle \Ds u_\eps,\f\rangle+o(1)
\le \langle (\Ds \psi)^+,\f\rangle+o(1),
\end{eqnarray*}
that readily gives $\Ds u\le (\Ds \psi)^+$ in the distributional sense in $\Omega$.

\bigskip

\noindent
{\bf Step 2} {\em Approximation argument.}

\medskip

Fix a small $\eps>0$ and put $\Omega_\eps:=\{x\in{\Omega}~|~\text{dist}(x,{\Omega})<\eps\}$. The convex set
$$
K_\eps=\{ v\in \widetilde H^s(\Omega_\eps)~|~v\ge \psi~~\text{a.e. on $\R^n$}~\}
$$
contains $K^s_\psi$, hence it is not empty. We denote by $u_\eps$ the unique solution
to the variational inequality
\begin{equation*}\tag{$\mathcal P_\eps$}
\label{eq:eps}
u_\eps\in K_\eps~,\qquad
\langle \Ds  u_\eps, v-u_\eps\rangle\ge 0\qquad\forall v\in K_\eps~\!,
\end{equation*}
so that $u_\eps\in \widetilde H^s(\Omega_\eps)$ and is nonnegative. Next we prove that
\begin{equation}
\label{eq:claim1}
0\le \Ds  u_\eps\le (\Ds \psi)^+\quad \text{in the distributional sense on $\Omega$.}
\end{equation}
For, we approximate $\psi$ in a standard way, via convolution. 
Let $(\rho_h)_h$ 
be a sequence of mollifiers such that $\text{supp}(\rho_h)\subset B_{\frac{1}{h}}$ and put
$\psi_{h}=\psi*\rho_h$. Notice that  for $h$ large enough,
$\psi_h= 0$ outside $\Omega_\eps$. Therefore 
\begin{equation}
\label{eq:psi_conv1}
\psi_h\in \widetilde H^s(\Omega_\eps)~,\qquad
\psi_h\to \psi\quad\text{in $ H^s(\R^n)$.}
\end{equation}
The convex set
$K_{\eps,h}:=\{v\in \widetilde H^s(\Omega_\eps)~|~u\ge \psi_h\}$ is not empty, as it contains $\psi_h$.
The variational inequality
\begin{equation*}\tag{$\mathcal P_{\eps,h}$}
\label{eq:vih}
u_h\in K_{\eps,h}~,\qquad
\langle \Ds  u_h, v-u_h\rangle\ge 0 \qquad\forall v\in K_{\eps,h}~\!,
\end{equation*}
has a unique solution $u_{h}\in \widetilde H^s(\Omega_\eps)$. Theorem \ref{T:Hs2} readily gives that
$u_h\to u_\eps$  in $\widetilde H^s(\Omega_\eps)$. Since
$\Ds \psi_h \in L^p(\R^n)$ for any $p\ge 1$, then Step 1 applies. In particular
\begin{equation}
\label{eq:eps_h}
0\le  \Ds  u_h\le (\Ds \psi_h)^+\quad\text{in the distributional sense on $\Omega$.}
\end{equation}
Next, $(\Ds \psi)^+*\rho_h$ is a nonnegative smooth function, and
$$(\Ds \psi)^+*\rho_h\ge (\Ds \psi)*\rho_h=\Ds\psi_h~\!.$$
Thus
$(\Ds \psi)^+*\rho_h\ge (\Ds\psi_h)^+$, and (\ref{eq:eps_h}) implies
$$
0\le  \Ds  u_h\le (\Ds \psi)^+*\rho_h\quad\text{in the distributional sense on $\Omega$.}
$$
Claim (\ref{eq:claim1}) follows, since $(\Ds \psi)^+*\rho_h\to (\Ds \psi)^+$ in the sense of measures,
and $\Ds u_h\to \Ds u_\eps$ in the sense of distributions. 

\bigskip

\noindent
{\bf Step 3} {\em Conclusion of the proof.}

\medskip

The last step in the proof consists in passing to the limit along a sequence $\eps\to 0$. First, we notice that
$u\in \widetilde H^s(\Omega_\eps)$ and in particular $u\in K_\eps$. 
Therefore, using the variational characterization of the unique solution $u_\eps$ to $(\mathcal P_\eps)$ we find
\begin{equation}
\label{eq:uff}
\frac{1}{2}\langle  \Ds  u_\eps,u_\eps\rangle\le 
\frac{1}{2}\langle  \Ds  u,u\rangle~\!.
\end{equation}
Now we fix $\eps_0>0$. 
Thanks to  (\ref{eq:uff}), we get that the sequence $u_\eps$ is bounded in $\widetilde H^s(\Omega_{\eps_0})$,
and therefore 
we can assume that $u_\eps\to \tilde u$ weakly in $\widetilde H^s(\Omega_{\eps_0})$.
From (\ref{eq:uff}) we readily get
\begin{equation}
\label{eq:uff2}
\frac{1}{2}\langle  \Ds  \tilde u,\tilde u\rangle\le 
\frac{1}{2}\langle  \Ds  u,u\rangle.
\end{equation}
On the other hand, $u_\eps\to \tilde u$ almost everywhere. Hence 
$\tilde u\in \widetilde H^s(\Omega)$ and $\tilde u\ge \psi$ on ${\Omega}$,
that is, $\tilde u\in K^s_\psi$.
Using the characterization of $u$ as the unique solution to the minimization problem (\ref{eq:minimization0}),
from (\ref{eq:uff2}),  (\ref{eq:uff}) we get that $\tilde u=u$ and $u_\eps\to u$  in $\widetilde H^s(\Omega_{\eps_0})$.
In particular, $\langle \Ds u_\eps,\f\rangle\to \langle \Ds u,\f\rangle$ for any $\f\in C^\infty_0(\Omega)$.
Now, from (\ref{eq:claim1}) we know that $( \Ds \psi)^+-\Ds {u_\eps}$ is a nonnegative distribution on $\Omega$.
Thus $( \Ds \psi)^+-\Ds {u}$ is nonnegative as well, and (\ref{eq:ineq}) is proved.
\QED

\medskip

\noindent
{\bf Proof of Theorem \ref{T:regularity}}\\
Statements $i)$ and $ii)$ hold by 
Theorem \ref{T:measure}. It remains to prove the last claim.

It is not restrictive to assume 
$f\equiv 0$. Hence $u$ solves $\mathcal P(\psi,0)$, 
 $\Ds u\ge 0$ by Theorem \ref{T:sup}, and $u$ is nonnegative in $\Omega$, 
see Remark \ref{R:MP}. Actually $u$ is lower semicontinuous
and positive by the strong maximum principle, see for instance
\cite[Theorem 2.5]{IMS}. Thus $u\ge \psi^+$ and $\{u>\psi\}=\{u>\psi^+\}$.

Next we use $c)$ in Theorem
\ref{T:sup} with $v= \psi^+\in \widetilde H^s(\Omega)$, to get
$$
\langle\Ds u,u-\psi^+\rangle =0.
$$
Let $\Omega'$ be any domain compactly contained in $\Omega$. We
claim that
\begin{equation}
\label{eq:zero}
\int\limits_{\Omega'} {\Ds u}~\!\!\cdot (u-\psi^+)~\!dx =0~\!.
\end{equation}
Since $\Ds u~\!\!\cdot (u-\psi^+)$ is a measurable nonnegative  function then the integral
in (\ref{eq:zero}) is nonnegative. To prove the opposite inequality
we put $g_m=(u-\psi^+)\wedge m$, $m\ge 1$. Let $\f$ be any nonnegative cut off function, 
with $\f\in C^\infty_0(\Omega)$
and $\f\equiv 1$ on ${\Omega'}$. 
Since $\Ds u\ge 0$, $\Ds u\in L^1_{\rm loc}(\Omega)$, $u-\psi^+\ge \f g_m$ and $\f g_m\in L^\infty(\Omega)$
has compact support in $\Omega$, we have that
$$
0=\langle \Ds u,u-\psi^+\rangle\ge\langle \Ds u,\f g_m\rangle=
\int\limits_{\Omega} {\Ds u}\cdot(\f g_m)dx\ge 
\int\limits_{\Omega'} {\Ds u}\cdot g_m dx.
$$
Next, use the monotone convergence theorem to get 
$$
0\ge \lim_{m\to\infty}\int\limits_{\Omega'} {\Ds u}~\!\!\cdot g_m~\!dx=
\int\limits_{\Omega'} {\Ds u}~\!\!\cdot (u-\psi^+)~\!dx,
$$
that concludes the proof of (\ref{eq:zero}).

Now, since $\Omega'$ was arbitrarily chosen and $\Ds u~\!\!\cdot (u-\psi^+)\ge 0$, equality  (\ref{eq:zero}) implies that
$\Ds u~\!\!\cdot (u-\psi^+)=0$  a.e. in $\Omega$,
and $iii)$ is proved.
\QED

\begin{Remark}
Theorem \ref{T:regularity} holds with the same proof also in the local case $s=1$.
Notice that no regularity assumptions on $\Omega$ are needed, and
the cases $p=1, p=\infty$ are included as well.
\end{Remark}

\begin{Remark}
\label{R:regularity}
To obtain better regularity results for $u$, one can apply the regularity
theory for
$$
\Ds u=g\in L^p(\Omega)\quad\text{in $\Omega$}~,\qquad u\in \widetilde H^s(\Omega).
$$
In particular, if $p>\frac{n}{2s}$ and $\Omega$ is Lipschitz and satisfies the exterior ball condition, then
$u$ is H\"older continuous in $\Omega$. See for example 
\cite[Proposition 1.4]{RoS} and \cite[Proposition 1.1]{ROS2}.
\end{Remark}

\bigskip

\noindent{\bf Proof of Theorem \ref{T:2}}\\
As usual, we can assume $f=0$.
Fix a small $\eps>0$, 
%
and let $\psi^\eps_{h}$ be a mollification of $\psi-\eps$. Then 
$\psi^\eps_{h}$ is smooth on $\overline\Omega$, 
$\psi^\eps_{h}<0$ on $\partial\Omega$ and 
$\psi^\eps_h\to \psi-\eps$ uniformly on $\overline\Omega$, as $h\to\infty$.

By Theorem \ref{T:regularity}, the solution $u_h\in \widetilde H^s(\Omega)$ to 
$\mathcal P(\psi^\eps_h,0)$ satisfies $\Ds u^\eps_h\in L^p(\Omega)$ and therefore 
$u^\eps_h$ is H\"older continuous, see Remark \ref{R:regularity}.
Moreover, the estimates in Theorem \ref{T:bounded1} imply
that $u^\eps_h \to u^\eps$ uniformly on  $\Omega$, where $u^\eps$ solves
$\mathcal P(\psi-\eps,0)$. In particular, $u^\eps\in C^0(\overline\Omega)$.
Finally, use again Theorem \ref{T:bounded1} to get that $u^\eps\to u$ uniformly, 
where $u$ solves $\mathcal P(\psi,0)$. In particular, $u$ is continuous on $\R^n$.

To check the last statement notice that the set $\{u>\psi\}\subseteq\Omega$ is open;
for any test function $\f\in C^\infty(\{u>\psi\})$ we have that $u\pm t\f\in K^s_\psi$ and therefore
$t\langle\Ds u,\pm \f\rangle\ge 0$ for
$|t|$ small enough. The conclusion is immediate.
\QED

\medskip
\noindent
{\bf Acknowledgments}.
R. Musina wishes to thank 
the National Program on Differential equations (DST, Government of India) 
and IIT Delhi for supporting her visit in January, 2015. A.I. Nazarov
is grateful to SISSA (Trieste) for the hospitality in October, 2015.

\end{document}